\documentclass[12pt]{article}
\usepackage[utf8]{inputenc}
\usepackage{amsmath}
\usepackage{amsfonts}
\usepackage{amssymb}
\usepackage{epsfig}
\usepackage{epstopdf}
\usepackage{indentfirst}
\usepackage{geometry}
\usepackage{graphicx}
\usepackage{caption}
\usepackage{subfig}
\usepackage[countmax]{subfloat}
\usepackage{rotating}
\usepackage{hyperref}
\usepackage{lscape}
\usepackage[numbers,sort&compress]{natbib}
\usepackage{tabularx}
\usepackage{setspace}
\usepackage{enumerate}
\usepackage{lineno}
\usepackage[section]{placeins}

\newtheorem{theorem}{Theorem}[section]
\newtheorem{lemma}{Lemma}[section]
\newtheorem{corollary}{Corollary}[section]

\newtheorem{definition}{Definition}[section]
\newenvironment{proof}{\smallskip\noindent{\textbf{Proof}}}{$\Box$}
\numberwithin{equation}{section}
\newcommand{\bmu}{{\boldsymbol \mu}}
\newcommand{\bet}{{\boldsymbol \eta}}

\newcommand{\blambda}{{\boldsymbol \lambda}}

\def\bbR{{\mathbb R}}

\def\by{{\boldsymbol y}}
\def\bh{{\boldsymbol h}}

\def\bv{{\boldsymbol v}}
\def\bg{{\boldsymbol g}}
\def\bA{{\boldsymbol A}}
\def\bB{{\boldsymbol B}}
\def\bC{{\boldsymbol C}}

\def\bW{{\boldsymbol W}}
\def\bx{{\boldsymbol x}}
\def\bc{{\boldsymbol c}}
\def\bC{{\boldsymbol C}}
\def\bX{{\boldsymbol X}}

\def\bv{{\boldsymbol v}}

\def\bu{{\boldsymbol u}}

\def\bq{{\boldsymbol q}}
\def\bp{{\boldsymbol p}}

\def\btau{{\boldsymbol \tau}}
\def\bmu{{\boldsymbol \mu}}
\def\bxi{{\boldsymbol \xi}}
\def\bbP{{\mathbb P}}

\def\bbI{{\mathbb I}}

\def\cI{{\mathcal I}}

\def\cK{{\mathcal K}}

\def\cS{{\mathcal S}}

\graphicspath{%
    {converted_graphics/}
    {/}
}

\begin{document}

\title{A geometry in the set of solutions to \\ ill-posed linear problems with box constraints: Applications to probabilities on discrete sets.}
\author{Henryk Gzyl \\
Center for Finance, IESA, Caracas.\\
{\small henryk.gzyl@iesa.edu.ve}
} 

\date{}
 \maketitle

\baselineskip=1.5 \baselineskip \setlength{\textwidth}{6in}

\begin{abstract}
When there are no constraints upon the solutions of the equation $\bA\bxi=\by,$ where $\bA$ is a $K\times N-$matrix, $\bxi\in\bbR^N$ and $\by\in\bbR^K$ a given vector, the description of the
set of solutions as $\by$ varies in $\bbR^K$ is well known. But this is not so when the solutions are required to satisfy $\bxi\in\cK\prod_{i\leq j\leq N}[a_j,b_j],$ for finite $a_j\leq b_j: 1\leq j\leq N.$ Here we provide a description of the set of solutions as a surface in the constraint set, parameterized by the Lagrange multipliers that come up in a related optimization problem in which
$\bA\bxi=\by$ appears as a constraint. It is the dependence of the Lagrange multipliers on the data vector $\by$ that determines how the solution changes as the datum changes.
 
The geometry on the solutions is inherited from a Riemannian geometry on the set of constraints induced by the Hessian of an entropy of the Fermi-Dirac type which is the objective in the restatement of the optimization problem mentioned above. We prove that the set of solutions is contained in $\ker(\bA)^\perp$ in the metric defined as the Hessian of the entropy.
\end{abstract}

\textbf{ Keywords}Inverse problems, Ill-posed problems, Convex constraints, Entropy Minimization, Hessian-Riemann metrics. 

\begin{spacing}{0.1}
   \tableofcontents
\end{spacing}

\section{Introduction and preliminaries}
Consider the typical linear inverse problem consisting of solving 
\begin{equation}\label{prob0}
\bA\bxi=\by
\end{equation}
 where $\bA:\bbR^N\to \bbR^K$ stands for both the linear transformation and its representing matrix in the canonical basis, and $\by\in\bbR^K$ is a given vector. In may problems $K\not= N$ and when equal, the matrix might not be invertible, this makes the problem an ill-posed, linear inverse problem. When there are no constraints upon the solution and the geometry on $\bbR^N$ and $\bbR^K$ is derived from the usual Euclidean norms, the space of solution as $\by$ varies is well known. 

For example, when $\by$ is in the range of $\bA$ (that is when \eqref{prob0} is consistent), if $\bA_g$ is any generalized inverse of $\bA$ satisfying $\bA\bA_g\bA=\bA,$ then $\bxi_g=A_g\by$ is a least squares solution to \eqref{prob0}.\\
When $\by$ is arbitrary and $\bA\bA_g$ is Hermitian, then $\bA_g\by$ is a least squares solution to \eqref{prob0}. See \cite{He} for example for proofs of these assertions and further references. The essence of the result is that the class of solutions is a linear subspace presented as the range of a linear mapping.

In many applications requires constraints to be imposed upon the solution. For examples problems in optimal transportation, image reconstruction, reconstructing probabilities from marginals or determining probabilities given the expected values of a few random variables, requires the imposition of convex constraints of box-type upon the solution. Problem \eqref{prob0} now becomes:

\begin{equation}\label{prob1}
 \begin{aligned}
&\bA\bxi=\by\\
\mbox{With}\;\;\; \bxi\in&\cK=\prod_{j=1}^N [a_j,b_j] \subset\bbR^N,\;\;\;\mbox{and}\;\;\by\in\bbR^K.
\end{aligned}
\end{equation}

The meaning of the symbols being as above, and the specification of the constraint set $\cK=\prod_{j=1}^N [a_j,b_j]$ depends on the specific application.To solve problem \eqref{prob1} we transform it into the following variational problem:

\begin{equation}\label{var1}
\mbox{Determine the point}\;\;\bxi^*=\inf\{\Psi(\bx)|\bxi\in \cK; \bA\bxi=\by\},
\end{equation}
where $\Psi$ is some appropriately specified convex function.\\
The remainder of the paper is organized as follows. After a rapid review of part of the large body of literature on the subject, in Section 2 we describe how to obtain the objective function $\Psi(\bxi)$ as the Fenchel-Lagrange dual of  the logarithm of the Laplace transform, denoted by $M(\btau),$ of a collection of unit masses at the vertices of the hypercube $\cK.$ We shall see how the product nature of the constraint set results in the objective function being an entropy of the Fermi-Dirac type for the objective function.  Besides of their roles in the solution of problem \eqref{var1}, the Hessians of $\Psi(\bxi)$ and $M(\btau)$ determine two related Riemannian structures, whose geodesics can be transported onto each other. Much of this material in convex analysis is known to information theorists, but is recalled here for the sake of completeness. 

In Section 3 we establish the representation of the solution to problems \eqref{prob1} and its variational version \eqref{var1}.  This representation serves to think of the set of solutions to \eqref{prob1}-\eqref{var1} as a surface in $\bbR^N.$ There we use the Bregman divergence defined by $\Psi$ as a heuristic aid to explain the nature of the solution to \eqref{var1}, and use the Bregman divergence to obtain some comparison results in Section 6. 

In Section 4  we finally consider the dependence of the solutions to \eqref{var1} in terms of the Lagrange multipliers. This allows us to think of the set of solutions to the inverse problem when $\by$ varies in the range of $\by$ as a hypersurface in $\bbR^N,$ which furthermore in the interior of $\bA(\cK).$ There we relate the Riemannian metric induced by the surface mappings to the Riemannian metric obtained by transporting by means of $\bA,$ the Hessian metric defined on the interior of $\cK$ by the Hessian of $\Psi.$ We also show that in this metric, the kernel $\ker(\bA)$ is orthogonal to the set of solutions to \eqref{prob1} as $\by$ varies in the interior of $\bA(\cK).$
There we also include the description of how do the Lagrange multipliers change when the datum changes.

In Section 5 we consider the particular case $\cK\subset[0,1]^N$ and some versions of the generalized moment problem in probability, that is, on how to determine a probability on a finite set if the expected value of a few random variables is known. We consider three variations on the theme: The first, besides the constraints imposed by the knowledge of the only constraints are the 
expected values of the observables. In the second example, we assume that previous knowledge of the probabilities is available in the form of an interval about an estimate. In the third example we suppose that the unknown probability has to be reconstructed from its marginals. This is similar, except for the normalization, to determine a supply and demand policy, when the supply and demand at a collection of sites is known. It is the first step of the transportation problem in which a cost function is to be minimized.

Section 6 contains some comparison results in which geodesic distances are compared to the standard Euclidean distance in the ambient space. We close this note with some final remarks.

For a recent application of this setup and maximum entropy in the mean (see next paragraph)) to the problem of reconstruction from marginals, plane tomography, the earth mover problem and supply and demand problems see \cite{Gz}. The same problem comes up in statistics when performing linear regression procedures with box constraints upon the unknown regression coefficients.

The procedure that we propose is related to the method of maximum entropy in the mean as described in \cite{DG} or \cite{GG}. It was actually used to provide a solution to the problem of reconstructing a contingency table from the knowledge of its marginals.  The main difference between the approach to be developed below and the method of maximum entropy in the mean, is that this uses the standard method of maximum entropy as an intermediate step, to determine an auxiliary probability , say $\bbP,$ supported by $\cK,$ which is the used to obtain the solution to \eqref{prob1} results from averaging a $\cK$ valued random variable with respect to $\bbP.$ Me mention that still another approach to transportation problem in which an entropy function is used as a regularizer in a linear programming problem is developed in \cite{Cu}. In \cite{Metal} a concept called the Tsallis entropy is used as a regularizer for the linear programming problem. 

To cite a few references intertwining the transportation problem, reconstruction of joint probabilities with given marginals and a geometry in the set of probabilities, using different approaches, we mention \cite{AC}, \cite{AKO}, \cite{Wong1}, \cite{WT} and \cite{KZ}.

\section{Geometries on $\bbR^N$ and on $\cK$}
\subsection{Mapping $\bbR^N$ onto $\cK$ and the objective function $\Psi(\bxi)$} 
As mentioned above, the starting point of our approach is the Laplace transform of a measure assigning a unit mass to each corner of the box $\cK.$ For that, write

\begin{equation}\label{ref}
dQ(\bxi) = \otimes_{j=1}^N\big(\epsilon_{a_j}(d\xi_j) + \epsilon_{a_j}(d\xi_j)\big).
\end{equation}

The notation $\epsilon_{a}(d\xi)$ stands for the point mass (Dirac) measure at the point $a.$ This reference measure is chosen for several reasons: First, its joint Laplace transform is trivial to compute and it plays a central role in the application of MEM to the inverse problem. Second, because its convex support is $\cK.$ For comparison results to the case in which $\cK$ is an arbitrary polytope, let $\bC=\{(c_1,...,c_{2N}): c_j\in \{a_j,b_j\}, j=1,...,N\}$ denote the set of corners of $\cK.$ Let us write $\delta_{\bc}(\bxi)=\prod_{j}\delta_{c_j}(dxi_j).$ Then \eqref{ref} becomes

\begin{equation}\label{ref1}
dQ(\bxi) = \sum_{\bc\in\bC}\delta_{\bc}(d\xi).
\end{equation}

With these notations, the Laplace transform of $dQ$, when the expression \eqref{ref} is used, is defined to be

\begin{equation}\label{LT1}
\zeta(\btau) = \int e^{\langle\btau,\bxi\rangle}dQ(\bxi) = \prod_{j=1}^N\big(e^{a_j\tau_j}+e^{b_j\tau_j}\big),\;\;\;\btau\in\bbR^N.
\end{equation}
Whereas, if \eqref{ref1} is used we obtain

\begin{equation}\label{LT2}
\zeta(\btau) = \int e^{\langle\btau,\bxi\rangle}dQ(\bxi)= \sum_{\bc\in\bC}e^{\langle\tau,\bc\rangle},\;\;\;\btau\in\bbR^N.
\end{equation}

Let us introduce the notation
\begin{equation}\label{mom1}
M(\btau) = \ln(\zeta(\btau)) = \sum_{j=1}^N\ln\big(e^{a_j\tau_j} + e^{b_j\tau_j}\big) \equiv \sum_{j=1}^N m(\tau_j).
\end{equation}
and call $M(\btau)$ the moment generating function (even though it is not the Laplace transform of a probability).  

The following lemma is used to obtain the Legendre-Fenchel dual of $M(\btau)=\ln\zeta(\btau).$ Since actually the variables $\tau_j$ appear separated in $M(\tau)$ we shall temporarily drop reference to the $j-$th subscript and put it back in when needed.

\begin{lemma}\label{lem1}
The function $\zeta(\tau):\bbR \to (a,b)$ given by $m(\tau)=\ln\big(\zeta(\tau)\big)$ is strictly convex and infinitely differentiable. Let $\xi\in (a,b),$ then the equation $\xi= dm(\tau)/d\tau$ has a unique solution that establishes a bijection between $\bbR$ and $(a,b).$
\end{lemma}
The proof is computational. Note that
$$\xi = \frac{\partial M(\btau)}{\partial \tau} = a\frac{e^{a\tau}}{e^{a\tau}+e^{b\tau}}+b\frac{e^{b\tau}}{e^{a\tau}+e^{b\tau}}$$
has a solution given by
\begin{equation}\label{chv}
e^{\tau} = \bigg(\frac{\xi-a}{b-\xi}\bigg)^{1/D}\;\;\;\Leftrightarrow\;\;\;\tau=\frac{1}{D}\ln\bigg(\frac{\xi-a}{b-\xi}\bigg).\;\;\; D=b-a.
\end{equation}
Note as well that $\xi(\tau) \to b$ (respectively to $a$) whenever $\tau\to-\infty$) (respectively to $+\infty$).
  
To sum up, the mappings 
$$\xi(\tau) = a\frac{e^{a\tau}}{e^{a\tau}+e^{b\tau}}+b\frac{e^{b\tau}}{e^{a\tau}+e^{b\tau}}\;\;\;\mbox{and}\;\;\;\tau=\frac{1}{D}\ln\bigg(\frac{\xi-a}{b-\xi}\bigg).\;\;\; D=b-a$$
are inverse to each other and are clearly infinitely differentiable. With this we put the subscripts back in (but we shall nevertheless drop them when it is redundant) and state:

\begin{definition}\label{chv}
Let us write $\Omega=int(\cK)=\prod_{1\leq j\leq N}(a_j,b_j).$ Define the mappings $\phi:\bbR^N \to \Omega$ and $\chi:\Omega \to \bbR^N$ whose components are specified by:
\begin{equation}\label{chv1}
\phi_j(\btau) = a_j\frac{e^{a_j\tau_j}}{e^{a_j\tau_j}+e^{b_j\tau_j}}+b_j\frac{e^{b_j\tau_j}}{e^{a_j\tau_j}+e^{b_j\tau_j}}\;\;\;\mbox{and}\;\;\;\chi_j(\bxi)=\frac{1}{D_j}\ln\bigg(\frac{\xi_j-a_j}{b_j-\xi_j}\bigg).\;\;\; D_j=b_j-a_j
\end{equation}
for $1\leq j \leq N.$
That is the mappings $\phi$ and $\chi$ are compositional inverses to each other.
\end{definition}

These changes of variable play a role in the computation of the Legendre-Fenchel dual of $M(\tau)$

\begin{lemma}\label{lem2}
Let $m(\tau)$ be defined as in \eqref{mom1}, that is $m(\tau)=\ln\big(e^{a\tau}+e^{b\tau}\big).$ Its Legendre-Fenchel dual is defined to be $\psi(\xi)=\sup\{\xi\tau-m(\tau): \tau\in\bbR\}$ for $\xi\in(a,b).$ It is given by:
\begin{equation}\label{LF1}
\psi(\xi) = \frac{\xi-a}{D}\ln\big(\frac{\xi-a}{D}\big) + \frac{b-\xi}{D}\ln\big(\frac{b-\xi}{D}\big).
\end{equation}
Finally, the Legendre-Fenchel dual of $M(\tau)$ is obtained by putting the subscripts back in, and it is given by:

\begin{equation}\label{LF2}
\Psi(\bxi) = \sum_{j=1}^N\frac{\xi_j-a_j}{D_j}\ln\big(\frac{\xi_j-a_j}{D_j}\big) + \frac{b_j-\xi_j}{D_j}\ln\big(\frac{b_j-\xi_j}{D_j}\big).
\end{equation}
\end{lemma}
Technically speaking, the function $\psi$ is defined as assuming the value $+\infty$ off the set $\Omega,$ but we shall not need that precision in this work.

The proof uses the results in Lemma \ref{lem1} and simple arithmetics to obtain (\ref{LF1}), and putting the subscripts back in, we obtain (\ref{LF2}). As usual it is consistent to take $0\ln 0=0.$ Some important, and simple to verify, properties of $\Psi(\bxi)$ are contained in the following result.

\begin{theorem}\label{propLF2}
The function $\Psi:\Omega\to\bbR$ defined in \eqref{LF2} is strictly convex, infinitely differentiable on $\Omega,$ continuous up to $\cK.$ Furthermore it achieves its maximum at the boundary, where it equals $0,$ and its minimum at the center of the box, at which it assumes the values $N\ln(1/2).$\\
Not only that, we also have the $N$-dimensional version of the remark after \eqref{chv1},
\begin{equation}\label{inv1}
\bxi = \nabla M_\tau(\btau) \Leftrightarrow \btau = \nabla_\bxi\Psi(\bxi).
\end{equation}
 \end{theorem}
\begin{corollary}\label{ortho0}
Let $H_M(\tau)$ and $H_\Psi(\bxi)$ denote the Hessian matrices of $M$ and $\Psi$ respectively. Then
\begin{equation}\label{ortho0.1}
H_M(\tau)H_\Psi(\bxi) = \bbI_N \;\;\;\mbox{whenever}\;\;\; \bxi = \nabla M_\tau(\btau).
\end{equation}
\end{corollary}

To finish this preamble, we determine the Bregman divergence associated to the convex function $\Psi(\bxi)$ and recall some of its basic properties to be used below.

\begin{theorem}\label{breg0}
The Bregman divergence determined by $\Psi$ is defined by
$$\delta_\Psi^2(\bxi,\bet) = \Psi(\bxi)-\psi(\bet)-\langle(\bxi-\bet,\nabla\Psi(\bet)\rangle.$$
For $\Psi(\bxi)$ as in (\ref{LF2}) we have
\begin{equation}\label{breg1}
\delta_\Psi^2(\bxi,\bet)= \sum_{j=1}^N\left[\frac{\bxi_j-a_j}{D_j}\ln\bigg(\frac{\bxi_j-a_j}{\eta_j-a_j}\bigg)+\frac{b_j-\bxi_j}{D_j}\ln\bigg(\frac{b_j-\bxi_j}{b_j-\eta_j}\bigg)\right].
\end{equation}
 Furthermore
\begin{equation}\label{propbreg}
\delta_\Psi^2(\bxi,\bet) \geq 0 \;\;\;\mbox{and} \;\;\;\delta_\Psi^2(\bxi,\bet)=0 \Leftrightarrow \bxi=\bet.
\end{equation}
\end{theorem}
$\Psi(\bxi)$ is reminiscent of an entropy function and $\delta_\Psi^2(\bxi,\bet)$ as a cross-entropy (or a Kullback divergence) defined on $\prod[a_j,b_j]$ if we think of $p(\xi)=(x-a)/(B-a)$ as the probability of picking a point uniformly on $(a,b).$ As usual, we

\subsection{The Riemannian structure on $\bbR^N$ induced by $M(\btau)$}
Again, since the coordinates are separated, we drop the subscripts and eventually put them back in.
The Hessian of $m(\tau)=\ln\big(e^{a\tau}+e^{b\tau}\big)$ is:
\begin{equation}\label{hess1}
h(\tau) = \frac{d^2 m(\tau)}{d\tau^2} = \left(\frac{b-a}{e^{d/2}+ e^{-d/2}}\right),\;\;\;d=\frac{b-a}{2}.
\end{equation}
To compute the distance between two points $\tau_0$ and $\tau_1,$ one needs to minimize the integral 

$$d_h(\tau_0),\tau_1) = \int_0^1 \bigg(h(\tau(t))(\dot{\tau})^2\bigg)^{1/2}dt$$

over all continuously differentiable curves $\tau(t)$ such that $\tau(0)=\tau_0$ and $\tau(1)=\tau_1.$
We will use the notation $\dot(\tau)$ as shorthand for time derivatives and, say $v'(tau)$ to denote derivatives with respect to the ``spatial'' variable. In the next section we show how to integrate the Euler-Lagrange the determines the geodesic. Here we just notice that if we write $h(\tau)= v'(\tau)^2$  we can write $h(\tau(t))(\dot{\tau})^2=\big(dv(\tau(t))/dt\big)^2.$ The function $v(\tau)$ that work here clearly is 

\begin{equation}\label{aux1}
v(\tau) = 2\arctan(e^{d\tau})
\end{equation}
This is clear since
$$v'(\tau) = 2d\frac{e^{d\tau}}{1 + e^{2d\tau}} = \frac{b-a}{e^{-d\tau} + e^{d\tau}} = h^{1/2}(\tau).$$
Since $h(\tau(t))(\dot{\tau})^2=\big(dv(\tau(t))/dt\big)^2,$ it is clear that the geodesics in the Riemannian metric associated to $h(\tau)$ are given by

\begin{equation}\label{geod1}
\tau(t) = v^{-1}\bigg(v(\tau_0) + t\big(v(\tau_1)-v(\tau_0)\big)\bigg).
\end{equation}
In this metric, the distance between two points $\tau_0$ and $\tau_1$ is:
$$d_h^2(\tau_0,\tau_1) = \big(v(\tau_1)-v(\tau_0)\big)^2.$$

To go back to the $N$-dimensional case, we write $\bh=H_M$ for the Hessian of $M,$ and think of it as of a (covariant) tensor field, that assigns the matrix $\bh\,$ (linear operator on the tangent space to $\bbR^N$ at $\btau\in\bbR^N$), with components 
\begin{equation}\label{hess1.1}
\bh_{n,m}(\btau)=h(\tau_n)\delta_{n,m}.
\end{equation}
We extend the $v$ introduced in \eqref{aux1} to a $\bv(\btau)$ which will be thought of as  vector as a vector valued function $\bv:\bbR^N\to \bbR^N$ given by $\bv(\btau)_n=v(\tau_n)$ whose Jacobian matrix, written $\bv'(\btau)$ for short, turn out to be diagonal and be the square root of $\bh,$ that is $\bh(\tau)=\bv'(\tau)^2.$

If we put the coordinate labels back in place, since the Riemannian metric is diagonal, the distance is the minimum of
$$d_{\bh}(\btau_0),\btau_1) = \int_0^1 \bigg(\sum_{j=1}^N h(\tau_j(t))(\dot{\tau_j})^2\bigg)^{1/2}dt$$
over all continuously differentiable curves $\btau(t)$ in $\bbR^N$ having $\btau(0)$ and $\btau(1)$ as initial and final points. Clearly, in this case \eqref{geod1} becomes:

\begin{equation}\label{geo2}
d_{\bh}(\btau(0)),\btau(0)) = \bigg(\sum_{j=1}^N \big(v(\tau_j(1))-v(\tau_j(0)\big)^2\bigg)^{1/2}.
\end{equation}
This is just the size of the velocity vector of the geodesic that goes from $\btau(0)$ to $\btau(1)$ in a unit of time. In general we have:

\begin{theorem}\label{int1}
With the notations introduced above, let $bX$ be a tangent vector to $\bbR^N$ at $\btau.$ Then
\begin{equation}\label{geod1.1}
\btau(t) = \bv^{-1}\bigg(\bv(\tau) + t\bX\bigg)
\end{equation}
is the geodesic curve starting at $\btau$ with velocity $\bX.$
Also, since in this case the geodesics are defined globally, the mapping
$$\exp_{\btau}:\bbR^N \to \bbR^N,\,\, \exp_{\btau}(\bX) = \bv^{-1}\bigg(\bv(\btau) + \bX\bigg)$$
is the geodesic map at $\btau.$
\end{theorem}
For a the basics about there concepts, see \cite{Hi} or \cite{S}.

\subsection{The Riemannian structure on $\Omega$ induced by $\Psi(\bxi)$}
The routine is as in the previous section. First note that, since the coordinates are separated in $\Psi(\bxi)$, drop the subscripts and consider the one dimensional case. To define a metric on $(a,b)$ by the Hessian of $\Psi,$ differentiating twice we have

\begin{equation}\label{hess2}
g(\xi) = \frac{d^2}{d\xi^2} \bigg(\frac{x-a}{D}\ln\frac{x-a}{D} + \frac{b-x}{D}\ln\frac{b-x}{D}\bigg) = \frac{1}{(\xi-a)(b-\xi)}.
\end{equation}

As above, put:
\begin{equation}\label{hess2.1}
g(\xi) = \big(u'(\xi)\big)^2\;\;\;\;\mbox{therefore}\;\;\;u'(\xi)=\bigg(\frac{1}{(\xi-a)(b-\xi)}\bigg)^{1/2}.
\end{equation}

Here, the computation specifying the $u$ is:

\begin{equation}\label{aux2.1}
u(\xi) = \int_a^{\xi}\frac{dx}{\sqrt{(x-a)(b-x)}} = \arcsin\left(\frac{2(\xi-m)}{D}\right)+\frac{\pi}{2}.
\end{equation}

Note that $u:(a,b)\to(0,\pi)$ smoothly. The change of variables $\xi\to\xi-m$ changes the denominator under the integral sign into 
$$\frac{D}{2}\sqrt{1-\big(\frac{2(\xi-m)}{D}\big)^2}$$
which reduces the integrand to a standard form. We put $m=(a+b)/2$ and recall that $D=b-a.$ 

As above, we follow the conventional notation $u'(\xi), u''(\xi)$ for derivatives with respect to ``spatial'' variables and $\dot{x}(t), \ddot{x}(t)$ for derivatives with respect to time variables. But this time we follow the traditional steps, we solve the $1-$dimensional case the Euler-Lagrange equations that determine the geodesic that goes from $\xi(0)$ at $t=0$ and $\xi(1)$ at $t=1,$ that is, we solve:

\begin{equation}\label{geo1}
\frac{d}{dt}\bigg(g(\xi)\dot{\xi}(t))\bigg)=\frac{1}{2}g'(\xi)\big(\dot{\xi}(t)\big)^2\;\;\;\Leftrightarrow\;\;\;g(\xi)\ddot{\xi}(t) + \frac{1}{2}g'(\xi)\big(\dot{\xi}(t)\big)^2 =0.\;\;\;
\end{equation}

Therefore 
$$\frac{d}{dt}\ln\bigg(\dot{\xi}g^{1/2}(\xi(t))\bigg)\;\;\; \Rightarrow u'(\xi)d\xi = \kappa dt\;\;\;\Rightarrow \xi(t)=u^{-1}\bigg(\xi(0)+\kappa t\bigg)$$
Above, $\kappa$ is an integration constant which can be determined from the initial and final data to be $\kappa=u(\xi(1))-u(\xi(0)).$ Therefore, on any interval $(a,b),$ the geodesic form $\xi(0)$ to $\xi(1)$ in the Riemannian metric \eqref{hess2.1}, is given by:

\begin{equation}\label{geo1.1}
\xi(t) = u^{-1}\big(u(\xi(0))+t(u(\xi(1))-u(\xi(0))\big), \;\;\;0\leq t \leq 1.
\end{equation}

We now pass onto the general case. We want to determine the geodesics in $\Omega$ when we use the Hessian $\bg =H_\Psi$ to define a Riemannian metric. Again, 
$\bg$ is diagonal, and invoking \eqref{hess2} we have:

\begin{equation}\label{hess2.2}
\bg_{i,j} = g(\xi_j)\delta_{i,j}.
\end{equation}

We introduce the vector valued function $\bu(\bx):\Omega\to (0,\pi)^N$ given by $\bu(\bxi)_j=u(\xi_j)$ with $u$ as in \eqref{aux2.1}. Denoting the Jacobian matrix of $\bv(\bxi))$ as a (diagonal) matrix $\bu'(\bxi)$, we write $\bg(\bxi)=(\bu'(\bxi))^2,$ and as in the previous section, we obtain:

\begin{theorem}\label{theo1}
Let $\bg$ defined in \eqref{hess2.2} be the Riemannian metric on $\Omega.$ 
The geodesic joining any two points $\bxi(1)$ and $\bxi(2)$ in $\Omega,$ is given by (the multidimensional extension of \eqref{geo1.1}) namely:

\begin{equation}\label{geo1.2}
\bxi(t) = \bu^{-1}\bigg(\bu(\bxi(0))+t(\bu(\bxi(1))-\bu(\bxi(0))\bigg), \;\;\;0\leq t \leq 1.
\end{equation}
The geodesic distance between any to points $\bxi(1)$ and $\bxi(2)$ in $\Omega$ is:
\begin{equation}\label{dist1}
d_{\bg}(\bxi(1),\bxi(2)) = \bigg(\sum_{j=1}^N\big(u_j(\xi_j(2)) - u_j(\xi_j(1))\big)^2\bigg)^{1/2}
\end{equation}
\end{theorem}

The result is clear since
$$\int_0^1\sqrt{\sum_{j=1}^N \big(u'(\xi_j(t))\dot{\xi}_j(t)\big)^2}dt = \int_0^1\sqrt{\sum_{j=1}^N \big(\frac{du(\xi_j(t))}{dt}\big)^2}dt = \bigg(\sum_{j=1}^N\big(u_j(\xi_j(2)) - u_j(\bxi_j(1))\big)^2\bigg)^{1/2}$$
where in the last step we used the fact that along the geodesic the velocity is constant and equals $\bv(\bxi(2))-\bv(\bxi(1)).$

\subsection{Transporting the geodesics}
Here we prove that the $\bh$-geodesics and the $\bg$-geodesics can be mapped onto each other. Let us begin with the following result.

\begin{lemma}\label{tg1}
Let the changes of variables $\Phi:\bbR^N\to \Omega$ and $\chi:\Omega\to\bbR^n$ be as in
 \eqref{chv1}, and let $\bv:\bbR^N\to\bbR^N$ and $\bu:\Omega\to\Omega$ be as in \eqref{aux1} and \eqref{hess2.1}. Then
\begin{equation}\label{tg2}
\bu(\phi(\btau)) = \bv(\btau)\;\;\;\;\mbox{and}\;\;\;\;\bv(\chi(\bxi)) = \bu(\bxi).
\end{equation}
\end{lemma}

The proof involves a simple computation that can be carried out in the one dimensional case due to the separability of the variables. For example, it is simple to verify that for any component:
$$u'(\Phi(\tau))\phi'(\tau) = v'(\tau) = \frac{(b - a)e^{\tau(a+b)/2}}{e^{a\tau}+e^{b\tau}}=\frac{(b - a)}{e^{-d\tau}+e^{d\tau}}.$$
The other case is similar. The result here is:

\begin{theorem}\label{tg3}
Let $\btau_0$ and $\btau_1$, respectively $\bxi_0=\phi(\btau_0)$ and $\bxi_1=\phi(\btau_1)$, be the initial and final points of $\bh$-geodesics, respectively, $\bg$-geodesics, in $\bbR^N$ and $\Omega,$ that are traversed in a unit of time. Then $\bxi(t)=\phi(\btau(t)).$
\end{theorem}
 
\begin{proof}
Again, it suffices to do the one-dimensional case. Consider the chain:
$$\xi(t) = u^{-1}\bigg(u(\xi_0) + t\big(u(\xi_1) - u(\xi_0)\big)\bigg)$$
Apply $u$ to both sides and use $\xi(t)=\phi(\tau(t))$ to obtain
$$u(\phi(\tau(t))) = u(\phi(\tau_0)) + t\big(u(\phi(\tau_1)) + t(u(\phi(\tau_1))-u(\phi(\tau_0))\big)$$
Now invoke Lemma \ref{tg1} to obtain
$$v(\tau(t)) = v(\tau_0) + t\big(v(\tau_1)-v(\tau_0)\big)$$
Hence:
$$\tau(t) = v^{-1}\bigg(v(\tau_0) + t\big(v(\tau_1)-v(\tau_0)\big)\bigg).$$
We leave it up to the reader to state and verify that the transport in the other direction can be verified similarly.
\end{proof}

\section{The solution the inverse problems \eqref{prob1}-\eqref{var1})}
At the beginning we argued on the passage form problem \eqref{prob1} to its variational version \eqref{var1}, that is to solve
\begin{equation}\label{prob3}
\mbox{Determine}\;\;\;\; \bxi^* \;\;\;\;\mbox{that realizes}\;\;\;\; \min\{\Psi(\bxi)| \bA\bxi=\by\}.\;\;\;\;\;\;\;\;\;\;\;\;\;
\end{equation}

Invoking Lagrange multipliers, the solution to (\ref{prob3}) is given by
\begin{equation}\label{repsol1}
\bxi_j^*= \frac{a_je^{a_j(\bA^t\blambda^*)_j} + b_je^{b_j(\bA^t\blambda^*)_j}}{e^{a_j(\bA^t\blambda^*)_j}+e^{a_j(\bA^t\blambda^*)_j}}.
\end{equation}
Here $\bA^t$ denotes the transpose of $\bA$ and $\blambda^*$ is a Lagrange multiplier that has yet to be determined. The heuristics of how to determine the Lagrange multiplier goes as follows. Use the property \eqref{propbreg} of the Bregman divergence introduced in \eqref{breg0}, with $\bet$ replaced by any point of the parametric surface  $\blambda\in\bbR^K\to\bxi(\blambda)$ in $\Omega$ given by:
$$  \bxi_j(\blambda) = \frac{a_je^{a_j(\bA^t\blambda)_j} + b_je^{b_j(\bA^t\blambda)_j}}{e^{a_j(\bA\blambda)_j}+e^{a_j(\bA^t\blambda)_j}}.$$
If $\bxi$ is any solution to \eqref{prob1} then

It takes a calculation to verify that if $\bxi$ is any possible solution to (\ref{prob3}),  then for any $\blambda\in\bbR^K$ we replace $\bet$ by $\xi(\blambda)$ in \eqref{breg1}, we obtain:
$$\delta_\Psi^2(\bxi,\bxi(\blambda)) \geq 0 \Leftrightarrow \Psi(\xi) \geq \langle\blambda,\by\rangle - M(\bA^t\blambda).$$

To finish (and in the next section we will make the statements more precise), note that the first order condition that determines the maximizer of $\Sigma(\by,\blambda)=\langle\blambda,\by\rangle - M(\bA^t\blambda)$ reads:
$$\by - \bA\nabla_{\tau}(\bA^t\blambda^*) = 0,$$
or, equivalently, $\bxi(\blambda^*)$ solves \eqref{var1}, and furthermore
\begin{equation}\label{duality}
\Psi(\bx(\blambda^*)) = \langle\blambda^*,\by\rangle - M(\bA^t\blambda^*).
\end{equation}

The interest in this argument lies in the fact that in order to solve (\ref{prob3}) it suffices to maximize a strictly convex function $\Sigma(\by,\blambda)$ over $\bbR^K.$

\section{The geometry on the class of solutions to \eqref{var1}}

To mimic the constructions in the previous section, we transport the measure $Q$ by means of $\bA$ onto $\bbR^K$ as usual. For a bounded, measurable function $f:\bbR^K\to\bbR$ we set: 

\begin{equation}\label{transpq}
\int f(\by)dQ_{\bA}(\by) = \int f(\bA\bxi)Q(d\bxi).
\end{equation}
Notice in passing that when $supp(f)\cap\bA(\bbR^N)=\emptyset$ then $int f(\by)Q_A(d\by)=0.$ Note that the Laplace transform of the transported measure $\zeta_{\bA}:\bbR^K\to\bbR$ is:

\begin{equation}\label{transLT}
\zeta_{\bA}(\blambda) = \int e^{\langle\by,\lambda\rangle}dQ_A(\by) = \int e^{\langle\bA\bxi\by,\lambda\rangle}dQ(\bxi) = \zeta(\bA^t\blambda).
\end{equation}

The comment right after \eqref{transpq} and the two exponentials in \eqref{transLT} suggest that there is no loss of generality in supposing that $\Omega\subset\ker(\bA)^{\perp}.$ We might as well suppose that the coordinate axes in $\bbR^K$ have been chosen so that the affine hull of $\bA(\cK)$ is contained in a subspace, and that $\blambda$ is parallel to that subspace. In other words, we might assume that, for the purposes of this approach, $\bA$ is onto $\bbR^K$ and that the relative interior of the convex set $\bA(\cK)$ is actually its interior. Let us now define the log-Laplace transform (moment generating function of $Q_{\bA})$ by

\begin{equation}\label{mom3}
M_{\bA}(\blambda) = \ln\big(\zeta_{\bA}(\lambda)\big) = M(\bA^t\blambda).
\end{equation}

\begin{lemma}\label{lem4}
With the notations just introduced, the function $M^{\bA}(\blambda)$ inherits the strict convexity and infinite differentiability from $M(\btau).$ Also, the Hessian matrix $\bh_{bA}$ of $M^{\bA}$ is
\begin{equation}\label{hess3}
\bh^{\bA}(\blambda) = \bA \bh(\bA^t\lambda)\bA^t.
\end{equation}
\end{lemma}

The proof is a computation left to the reader.\\
Let now $\blambda: \cI\to\bbR^K$ be a continuously differentiable curve defined on an interval $\cI$ containing $(0,1).$ Observe that making use of \eqref{hess3} the Lagrangian function yielding the geodesics of $\bh^{\bA}$ is 

\begin{equation}\label{lag1}
\langle\dot{\blambda},h^{\bA}(\blambda)\dot{\blambda}\rangle = \langle \bA^t\dot{\blambda}, h(\bA^t\blambda)\bA^t\dot{\blambda}\rangle.
\end{equation}
Now, recall from Section (2.2) that $\bh$ is diagonal with entries $h_{jj}(\btau)=v'(\tau_j))^2$ where $\bv$ is detailed in \eqref{aux1}. Thus, \eqref{lag1} can be rewritten as

\begin{equation}\label{lag2}
\langle\dot{\blambda},\bh^{\bA}(\blambda)\dot{\blambda}\rangle = \langle\frac{d}{dt}\bv(\bA^t\blambda),\frac{d}{dt}\bv(\bA^t\blambda)\rangle,\;\;\;\mbox{with}\;\;\;\bv_j\big((\bA^t\blambda)_j\big).
\end{equation} 

And we now state:
\begin{theorem}\label{geotr1}
Suppose that the curve $\blambda: \cI\to\bbR^K,$ is continuously differentiable on an interval $\cI$ containing $(0,1),$ is a geodesic in $\bbR^K$ in the metric $h_{\bA},$ passing through $\blambda_0$ and $\blambda_1$ at $t=0$ and $t=1$ respectively. Then the curve $\btau(t)=\bA^t\blambda(t)$ is a geodesic in the $d_{\bh}$ distance on $\bbR^N,$ that passes through $\btau_0=\bA^t\blambda_0$ at $t=0$ and through $\btau_1=\bA^t\blambda_1$ at $t=1.$
\end{theorem}

\begin{proof} 
After the change of coordinates $\bet=\bA^t\blambda$ the distance along the curve $\bet(t)$ is 
$$\int_0^1 \big(\langle\dot{\bxi},\dot{\bxi}\rangle\big)^{1/2}dt.$$
This is minimized by $\bxi(t)=\bxi(0)+t\big(\bet(1)-\bet(0)\big).$ Substituting the right symbols  we obtain
\begin{equation}\label{geotr1.2}
\bv(\bA^t\blambda(t)) = \bv(\bA^t\blambda(0)) +t\big(\bv(\bA^t\blambda(1)) - \bv(\bA^t\blambda(0))\big).
\end{equation}
Or
$$\bA^t\blambda(t)=\bv^{-1}\bigg(\bv(\bA^t\blambda(0)) + t\big(\bv(\bA^t\blambda(1)) - \bv(\bA^t\blambda(0))\big)\bigg).$$
According to \eqref{geo1}, $\bA^t\blambda(t)$ is a geodesic trough the given initial and final points.
\end{proof}

This result con be combined with Theorem \ref{tg3}, to obtain the geodesics in the space of solutions as follows.

\begin{theorem}\label{key}
As above, let $\blambda: \cI\to\bbR^K,$ be continuously differentiable on an interval $\cI$ containing $(0,1),$ is a geodesic in $\bbR^K$ in the metric $\bh_{\bA},$ passing through $\blambda_0$ and $\blambda_1$ at $t=0$ and $t=1$ respectively. Then $\bxi(\blambda(t))$ given by is a geodesic in $\cS$ and 

\begin{equation}\label{geosol}
\bxi(\blambda(t)) = \bu^{-1}\bigg(\bu(\bxi(\blambda(0))) + t\big(\bu(\bxi(\blambda(1)))-\bu(\bxi(\blambda(0)))\big)\bigg).
\end{equation}
\end{theorem}

The computation of $\Psi_{\bA}(y),$ the Fenchel-Legendre of $M_{\bA}(\lambda)$ has an interesting byproduct, namely, the solution to \eqref{var1}.

\begin{theorem}\label{dualLT2}
Let $M_{\bA}(\blambda)$ be as in \eqref{mom3}. Let $\by$ be a interior point of the closed convex  set $\bA(\cK).$Suppose that
\begin{equation}\label{LF2}
\Psi_{\bA}(\by) = \sup\{\langle\by,\blambda\rangle - M_{\bA}(\blambda):\blambda\in\bbR^K\} = \Psi(\bA_g\by).
\end{equation}
is bounded above and let $\blambda^*$ the point at which its maximum is reached and 
\begin{equation}\label{sol1}
\bxi^* = \nabla_{\btau}M(\blambda^*)
\end{equation}
is a solution to \eqref{var1}.
\end{theorem}

The interesting remark at this point is that the statement of the theorem contains the prescription to determine $\blambda^*$ numerically.

\begin{proof} 
That a maximizer exists follows from Theorem \ref{propLF2} and the fact that, since $\by=\bA\bxi_o$ for some $\bxi_o\in\Omega,$ then 
$$\Psi_{\bA}(\by) = \sup\{\langle\by,\blambda\rangle - M_{\bA}(\blambda):\blambda\in\bbR^K\} \leq \Psi(\bxi_o) \leq 0.$$
Observe that the condition for $\blambda^*$ to be a maximizer is that $\by-\bA\nabla_{\btau}M(\blambda^*)$ 
\end{proof}

So, turning this around, what \eqref{repsol1} or \eqref{sol1} assert, is that the mapping $phi:\bbR^N\to\Omega$ restricted to $\bA(\bbR^K) = (\ker(\bA)^\perp$ provides us with a description of the class solutions to $\bA\bxi=\by$ as $\by$ ranges in the interior of $\bA(\cK).$ So, the mapping:

\begin{equation}\label{surf1}
S:\bbR^K \to \Omega,\;\;\;\;\mbox{given by}\;\;\;S(\blambda)=\phi(\bA^t\blambda)= \bxi(\lambda)
\end{equation}
 describes a parameterized surface in $\bbR^N..$ Let us call that surface $\cS$ to give it a name. To define a metric on the tangent space to $\cS$ we start with:.

\begin{theorem}\label{surf2}
Let $\blambda\in\bbR^K \to \Omega$ be the map defined in \eqref{sol1}.  The metric on the tangent space to $\cS$ is defined by
$$G_{k,l} \equiv \langle\frac{\partial \bxi}{\partial\lambda_k},\frac{\partial \bxi}{\partial\lambda_l}\rangle_{\bg}$$\
Then
\begin{equation}\label{geo3}
G_{k,l} = \sum_{n.m}A_{k,n}\frac{\partial M}{\partial\tau_n\partial\tau_m}A_{l,m}.
\end{equation}
Or in matrix notation $G(\blambda)=\bA \bh(\bA^t\blambda)\bA^t$ where, recall, $h$ stands for the Hessian matrix of $M.$ 
\end{theorem}

Observe that this is consistent with the approach leading to \eqref{hess3}. The proof is a computation that uses the fact that the Hessian matrices $\bg$ and $\bh$ are (multiplicative) inverses of each other as mentioned in \eqref{ortho0.1}.

Comparing to \eqref{hess3} we see that there are two ways of putting a metric on $\bbR^K.$thinking about the metric $G.$ The one made explicit in Theorem \ref{surf2} is consistent with thinking on the space of solutions to \eqref{prob1} as a manifold in $\cK$ coordinatized by $\bbR^K.$

Furthermore we also have;

\begin{theorem}\label{ortho}
Let $\blambda\in\bbR^K \to\Omega$ be the map defined in \eqref{sol1}. Then $\cS\subset \ker(\bA)^{\perp}.$ That is, in the metric $\bg$ of the ambient space $\bbR^N,$ if $\bv\in\ker(\bA),$  $\bv$ be a vector field in $\bbR^N$ defined on $\cS,$ and let $t\to\blambda(t)$ be a differentiable curve in $\bbR^K.$ Then
 then 
$$\langle\bv,\frac{d\bxi(\blambda)}{dt}\rangle_{\bg} = 0.$$
\end{theorem}

\begin{proof}
The proof is computational. Let $\bv$ be a vector field in $\bbR^N$ defined on $\cS,$ and let $t\to \blambda(t)$ be a differentiable curve in $\bbR^K,$ and $\bxi(t)$ be its image in $\bbR^N$ by \eqref{surf1}. Then, since the matrices $\bg$ and $\bh$ are inverses to each other, we have:

$$\langle\bv,\frac{d\bxi}{dt}\rangle_{\bg} = \langle\bv, \bg\bh\bA^t\dot{\blambda}\rangle = \langle \bA\bv,\dot{\blambda}\rangle.$$
This certainly vanishes when $\bv\in\ker(\bA).$
\end{proof}

At this point we make use of the results in Section 2.4, Theorem \ref{geotr1} and \eqref{geotr1.2} to state a different version of Theorem \ref{key}:

\begin{theorem}\label{geotr3}
Let $\blambda(t)$ a geodesic in $\bbR^K$ as in Theorem \ref{geotr1}. Then $\xi(\blambda(t))$ is a geodesic lying on the hypersurface $\cS$ 
\end{theorem}

\begin{proof} Observe that, as above,
$$\langle\frac{d\bxi}{dt},\frac{d\bxi}{dt}\rangle_{\bg} = \langle\dot{\blambda}\bA,\bh(\bA^t\blambda)\bA^t\dot{\blambda}\rangle $$
where $\bh(\btau)$ is the Hessian matrix of $M(\btau).$ Using the factorization of $\bh$ as above, we see that 
$$\langle\frac{d\bxi}{dt},\frac{d\bxi}{dt}\rangle_{\bg} = \langle\frac{d\bv(\bA^t\blambda)}{dt},\frac{d\bv(\bA^t\blambda)}{dt}\rangle,$$
and we already saw that the right hand side is minimal when $\blambda(t)$ are geodesics in  $\bbR^K.$
\end{proof}

To state an interesting role of the metric $G,$ notice that $\blambda$ depends non-linearly on $\by$ through \eqref{sol1}. From that we obtain the following result.
\begin{lemma}\label{surf4}
Let $\bA\bxi(\blambda)=\by.$ Let $\delta\blambda/\delta\by$ denote the Jacobian matrix with components $\partial\lambda_n/\partial y_m$ for $1\leq n,m\leq K.$ then
\begin{equation}\label{surf5}
\frac{\delta \blambda}{\delta \by} = G^{-1}(\blambda).
\end{equation}
\end{lemma}
\begin{proof}
Let $\bxi(\blambda)$ be as in \eqref{sol1}. Differentiate both sides of $\bA\bxi(\lambda) = \by$ with respect to $\by$ and invoke \eqref{geo3} to obtain 
$$\bA \bh(\bA^t\blambda)\bA^t \frac{\delta \blambda}{\delta \by} = G(\blambda)\frac{\delta \blambda}{\delta \by}=\bbI.$$
Since the determinant of the identity  is $1,$ it implies that the determinant of each factor on the left hand side does not vanish, thus the claim follows.
\end{proof}

A potentially useful application of this result goes as follows:
\begin{corollary}\label{surf6}
With the notations introduced above, suppose that $\by$ changes as $\by\to\by+\Delta\by$ Then, up to first order in $\Delta\by,$ the solution to the inverse problem \eqref{prob1} changes as
\begin{equation}\label{sol2}
\Delta\bxi = \bh(\bA^t\blambda)\bA^tG^{-1}\Delta\by.
\end{equation}
\end{corollary}

\section{Applications to probabilities on discrete sets}
As mentioned in  the introduction we will consider a few variations on the theme of determining a probability on a set of $N$ points when the expected value of $K<N$ random variables is prescribed and that a prescribed range for the unknown probabilities is available. In all examples below, the name of the unknown variables changes from $\xi_j$ to $p_j.$

\subsection{Case 1: The Only the expected value of a few observables is known}
The problem to solve consists of

\begin{gather}
\mbox{Determine}\;\;\;\{p_1,...p_N\}\;\;\;\mbox{such that:} \label{ex1.1}\\ 
\sum_{j=1}^N p_j = 1\;\;\;\mbox{and}\;\;\; \label{ex1.2} \\
\sum_{j=1}^N B_{k,j}p_j = y_k,\;k=2,...,K.\label{ex1.3}
\end{gather}
For this example we choose $a_j=0$ and $b_j=1$ for $j=1,...,N.$ The matrix $\bA$ for this example is
\begin{equation}\label{mata}
\bA = {\bu^t\atopwithdelims[]\bB},\;\;\;\Rightarrow\;\;\;\bA^t=[\bu\,\bB]
\end{equation}
where $\bu$ denotes the $N$-vector with all entries equal to $1.$ We will use $\blambda^1=(lambda_1,\bmu^t)$ with $\bmu\in\bbR^{K-1}.$ Note as wee that in this (and all cases below), we have $\big(\bA^t\blambda\big)_k =\lambda_1+\big(\bB^t\bmu\big)_k.$
With these notations, it is clear that in this case:
$$M(\bA^t\blambda) = \sum_j \ln\big(1 + e^{\lambda_1+(\bB^t\bmu)_j}).$$
Therefore the probabilities that solve the problem are obtained differentiating $M(\bA^T\blambda)$ with respect to $\lambda_1$ and are given by :\begin{equation}\label{up1}
p_j = \frac{e^{\lambda_1+(\bB^t\bmu)_j}}{1 + e^{\lambda_1+(\bB^t\bmu)_j}}.
\end{equation}
which of course satisfy the constraints:
$$\sum_{j=1}^N B_{k,j}p_j = y_k.\;\;\;k=2,...,K.$$
The important thing is that for numerical purposes one does not need to solve a nonlinear set of equations to determine the Lagrange multipliers, but to use a gradient based minimization to find the minimizer of $M(\bA^T\blambda)- \lambda_1-\langle\by,\bmu\rangle.$

\subsection{Case 2: The probabilities are known to lie in a range}
This time the problem to solve is:

\begin{gather}
\mbox{Determine}\;\;\;\{p_1,...p_N\}\;\;\;\mbox{such that:} \label{ex2.1}\\ 
p_j \in (p^0_j-\delta_j,p^0_j+\delta_j),\; j=1,...,N, ;\;\mbox{and} \label{ex2.2}\\
\sum_{j=1}^N B_{k,j}p_j = y_k,\;k=2,...,K.\label{ex2.3}
\end{gather}
The $p^0_j, j=1,...,N$ are taken to be probabilities obtained in some previous experiment, or provided by some theoretical model, and clearly, the interval does not have to be symmetric about $p^0_j,$ but the notations are a trifle simpler. 
In this case the integral constraints consist of the expected vale of the a few observables, there is no constraint as in \eqref{ex1.2}, and the solution the problem is given by

\begin{equation}\label{up2}
p_j = \frac{(p^0_j-\delta_j)e^{(p^0_j-\delta_j)\bA^t\blambda)j} + (p^0_j+\delta_j)e^{(p^0_j+\delta_j)\bA^t\blambda)j}}{e^{(p^0_j-\delta_j)\bA^t\blambda)j} + e^{(p^0_j+\delta_j)\bA^t\blambda)j}}
\end{equation}
Again, th is to be obtained $\blambda$ by minimizing the appropriate strict convex function over $\bbR^K.$

\subsection{Case 3: Reconstruction of joint probabilities from marginals}
This problem is stated as follows.

\begin{gather}
\mbox{Determine}\;\;\;\{P(i,j)|\,1\leq N,\,1\leq M\}\;\;\;\mbox{such that:} \label{ex3.1}\\ 
\sum_{j=1}^M P(i,j)=P_i,\;i=1,...,N,\;\;\mbox{and}\;\;\;\sum_{i=1}^N P(i,j)=Q_j,\;j=1,...,M. \label{ex3.2}\\
\sum_{i=1,j=1}^{N,M} W(i,j)P(i,j)= w\label{ex3.3}
\end{gather}

Problem \eqref{ex3.1} subject to the constraints \eqref{ex3.2} is the standard problem consisting of reconstructing a joint probability from its marginals. Observe that since the sum of the marginal probabilities is $1,$ no further normalization constraint is necessary. When the constraint \eqref{ex3.3} is added, the problem becomes a cost constrained transportation problem. As we want to examine the geometry of  \eqref{ex3.1}-\eqref{ex3.2} first, let us do the vectorization of the problem in two stages. 

For that we list the elements of  $\{1,...,N\}\times\{1,...,N\}$ in lexicographic order. that is,to the joint probability $P(i,j)$ we associate a $NM$-vector $\bx\in \cK$ (or $\bxi\in [0,1]^{NM}$ if no constraints on the size of the probabilities are imposed). 
  
\begin{equation}\label{ex3.4}
\bC\bx = {\bp \atopwithdelims()\bq} = \by,\;\;\;\mbox{where}\;\;\;\bC={\bC_1 \atopwithdelims()\bC_2}.
\end{equation}
\noindent  Above, $\bC_1$ stands for an $N\times NM$-matrix taking care of the row constraints, $\bC_2$ denotes an $M\times NM$ matrix taking care of the column constraints, and  $\bp=(P_1,...,P_N)^t$ and $\bq=(Q_1,...,Q_M)^t,$ denote, respectively the row and column sums. As usual the superscript $t$ denotes transposition (we think of vectors as columns). How to transform \eqref{ex3.2} into the action of matrix upon $\bxi$ is explained in \cite{Gz}.

The linear cost constraint consists of the specification of an $NM$-vector $\bW$ and require the solution to problem  \eqref{ex3.1}-\eqref{ex3.2} to satisfy the constraint: 

\begin{equation}\label{ex3.5}
\sum_{n=1}^{NM} W_{n}\pi_{n} = w.
\end{equation}

Now, the augmented problem  \eqref{ex3.1}-\eqref{ex3.3} in vectorized form becomes : 

\begin{equation}\label{prob3.6}
\begin{aligned}
\mbox{Determine}\;\;\; \bxi^*\in&[0,1]^{NM}\;\;\;\mbox{such that:}\\
\bA\bxi^* &= \by(w)
\end{aligned}
\end{equation}

where explicitly: 

\begin{equation}\label{augm1}
 \bA ={\bC \atopwithdelims[]\bW^t},\;\;\;\by(w) = {\by \atopwithdelims[] w}.
\end{equation}

Also, besides prior information upon the size of the $\xi_j,$ the box constraints $\xi_j\in[a_j,b_j]$
may have a public policy explanation if the problem is interpreted as a supply demand problem. Not only that, one may use $\nu$ as a control parameter to obtain a solution with as small as possible.

For this example it is typographically to distinguish the Lagrange multipliers. We will write the multipliers as $(\blambda,\bmu,\nu)^t,$ where $\blambda\in\bbR^N,$ $\bmu\in\bbR^M$ and $\nu\in\bbR$ for the multipliers corresponding to the sum by row and sum by column constraints, and the cost constraint respectively. With these notations, the solution to problem \eqref{ex3.6} is:

\begin{equation}\label{ex3.7}
\bxi_j =  \frac{a_je^{a_j(\lambda_j+\mu_j+\nu)} + b_je^{b_j(\lambda_j+\mu_j+\nu)}}{e^{a_j(\lambda_j+\mu_j+\nu)} + e^{b_j(\lambda_j+\mu_j+\nu)}}
\end{equation}

 To find an optimal transport policy is  complete by a sequence of similar problems for a decreasing sequence of costs.

\section{Some comparison results}
Let us begin by the upper bound. As usual, we drop the subscripts until we need them. Notice that
$$\sqrt{(\xi-a)(b-\xi)} \geq \frac{1}{D}(\xi-a)(b-\xi).$$
Therefore
$$D\int_\xi^\eta \frac{dx}{(\xi-a)(b-\xi)} = \int_\xi^\eta u'(x)dx = u(\eta)-u(\xi).$$
But
$$D\int_\xi^\eta \frac{dx}{(\xi-a)(b-\xi)}=\int_\xi^\eta\bigg(\frac{dx}{x-a}+\frac{dx}{b-x}\bigg)=\ln\left(\frac{b-\xi}{\xi-a}\right) -\ln\left(\frac{b-\eta}{\eta-a}\right).$$
after rearranging the right hand side a bit. We gather these computations under
\begin{theorem}\label{comp1}
For $\bx,\bet\in\Omega,$ the computations made above imply:
\begin{equation}\label{comp1.1}
\bigg(\sum_{j=1}^N\left(\ln\big(\frac{b_j-\xi_j}{\xi_j-a_j}\big)-\ln\big(\frac{b_j-\eta_j}{\eta_j-a_j}\big)\right)^2 \bigg)^{1/2}  \geq d_{\bg}(\bxi,\bet).
\end{equation}
\end{theorem}
To establish a curious collection of lower bounds note the following. Centering at $m=(a+b)/2$  and putting $z=(2\xi-m)/D,$ we obtain:
$$u'(\xi) = \frac{2}{D}\sqrt{\frac{1}{(1 - z^2)}} \geq\left\{\begin{matrix}
                                                       \sqrt{1+z^2}\\
                                                        |z|^n\;\;\;n\geq 0.\end{matrix}\right.$$
																								
The right hand side follows from neglecting as many terms in the expansion $(1-z^2)^{-1}=\sum_{n\geq 0}z^{2n}$ as desired. Note as well that when $\xi(2)>\xi(1)\;\Rightarrow\;z(2)>z(1),$ we have
$$\int_{z(1)}^{z(2)}\sqrt{1+z^2}dz=\sinh(z(2)) - \sinh(z(1)).$$
Similarly, when $n=0$ we have:
$$\int_{z(1)}^{z(2)}dz = \xi(2) - \xi(1).$$
These two cases yield the following comparison result.
\begin{theorem}\label{comp2}
With the notational conventions just introduced, for $\bx,\bet\in\Omega$ we have:
\begin{gather}
d_{\bg}(\bx,\bet) \geq \sqrt{\sum_{j=1}^N\left(\frac{2}{D_j}\right)\left[\sinh\left(\frac{2(\eta_j-m_j)}{D_j}\right)
-\sinh\left(\frac{2(\eta_j-m_j)}{D_j}\right)\right]^2 }\label{comp2a} \\
d_{\bg}(\bx,\bet) \geq \sqrt{\sum_{j=1}^N\left(\eta_j - \xi_j\right)^2}.\label{comp2b}
\end{gather}
\end{theorem}

\section{Final remarks}
The results obtained above rely heavily on the fact that the constraint space is a box. Note that if we label the corners of the box as $\bC=\{\bc=(c_1,c_2,...,c_N)|c_j\in\{a_j,b_j\}\}$, and if we denote by $\delta_{\bc}(d\bxi)$ the unit mass measure supported by $\bc\in\bbR^N,$ since 
$$\delta_{\bc}(d\bxi) = \prod_{j=1}^N\delta_{c_1}(d\xi_1)(d\xi_2)...(d\xi_N)$$
then the generating function $M(\btau)$ can be written as
\begin{equation}\label{gen1}
M(\btau) = \ln\bigg(\sum_{\bc\in\bC}e^{\langle\bc,\btau\rangle}\bigg) = \ln\zeta(\btau).
\end{equation}

Up to this point, this is generic nd can be done for any polytope in $\bbR^K$ instead of the box $\cK.$ In the general case, one does not have the factorization that leads to \eqref{LT2} and to \eqref{mom1}. Nevertheless, when the polytope is the unit simplex in $\bbR^N,$ instead of the geometry developed in Section 2, one still has a possible geometry, one still has a geometry in the simplex inherited from a geometry on the strictly positive orthant in $\bbR^N$ considered as a commutative group. The details of this approach are explained in \cite{GN}. Here we only mention that in the resulting geometry, the geodesic joining the probability mass distributions $(p_1,....,p_N)$ and $(q_1,...,q_N)$ is given by
$$p_j(t) = \frac{p_j^{1-t}q_j^{t}}{\sum_{n}p_n^{1-t}q_n^{t}},\;\;\;0\leq t \leq 1.$$
To further stress the difference between the two cases, note that the hypercube need not be contained in the positive orthant of $\bbR^N.$ Note as well that in \eqref{gen1} the corners of the hypercube can be replaced by the corners (or the extreme points) of any convex polytope in $\bbR^N.$ The question is: which convex sets lead to computations like those carried out in Section 2?

\textbf{Compliance, ethical issues, conflicts of interest and data availability} This work was prepared without recurring to AI in any way, it involves no conflicts of interests, it was written only by the corresponding author, and it does not use any data.


\begin{thebibliography}{}
\bibitem{A} Amari, S-I. \textit{``Information Geometry and its Applications''}, Springer (2016).
\bibitem{AC} Amari,S-I and Cichocki, A. (2010). {\it Information geometry of divergence functions}, Bulletin of the Polish Academy of Sciences, \textbf{58}, 183-195.
\bibitem{AKO} Amari, S., Karakida, R. and Oizumi, M. (2018). {\it Information geometry connecting Wasserstein distance and Kullback-Leibler divergence via the entropy relaxed problem}, Information Geometry, \textbf{1}, 13-37.
\bibitem{Cu} Cuturi, M. (2013). {\it Sinkhorn distances: lightspeed computation of optimal transport}, Proceedings of the 26th International Conference on Neural Information Processing Systems, \textbf{2}, 2292-2300.
\bibitem{DG}Dacunha-Castelle, D. and Gamboa, F. (1990). {\it Maximum d'entropie et probleme des moments}, Annals de l'Institut Henri Poincar\'e, \textbf{26}, 567-596.
\bibitem{GG} Golan, A and Gzyl, H (2002). {\it A generalized maxentropic inversion procedure for noisy data},  Applied Math. and Computation, \textbf{127}, 249-260.
\bibitem{G1} Gzyl. H. (2019). {\it Best predictors in logarithmic distance between positive random variables}, Journal Applied Mathematics, Statistics and Information,  \textbf{15}, 15-28. 
\bibitem{GN} Gzyl, H. and Nielsen, F. (2020). {\it Geometry of the probability simplex and its connection to the maximum entropy method}, Journal of Applied Mathematics, Statistics and Information, \textbf{16}, 25-35.
\bibitem{Gz} Gzyl, H. (2020). {\it Construction of contingency tables by maximum entropy in the mean}, Communications  in Statistics: Theory and Methods, \url{https://doi:10.1080/03610926.2020.1723639}.
\bibitem{He} Hearon, J.Z. (1968). {\it Generalized inverses and solutions to linear system}, Journal of Research of the National Bureau of Standards-B, Mathematical Sciences, \textbf{178}, 303-308.
\bibitem{Hi} Hicks, N.J. (1965). {\it Notes of Differential Geometry}, D. Van Nostrand Company Inc., Princeton.
\bibitem{Jaynes} Jaynes, E. T. (1957). {\it Information theory and statistical mechanics}. Physical review, 106(4), 620.  
\bibitem{KZ} Khan, G., Zhang, J. (2022). {\it When optimal transport meets information geometry}, Information. Geometry \textbf{5}, 47-78.
\bibitem{Metal} Muzellec, B., Nock, R., Patrini, G, and Nielsen, F. (2016). {\it Tsallis Regularized Optimal Transport and Ecological Inference},  Proceedings of the Thirty-First AAAI Conference on Artificial Intelligence, 2387-2393. \url{DOI:10.1609/aaai.v31i1.10854}
\bibitem{S} Sternberg, S. (2012). {\it Curvature in Mathematics and Physics}, Diver Publications Inc., Mineola, N.Y.
\bibitem{Wong1} Wong, T-K.L. (2018). {\it Logarithmic divergences from optimal transform and R\'enyi geometry}. Information Geometry, \textbf{1}, 39-78.
\bibitem{WY} Wong, T-K.L. and Yang, J. (2022). {\it Pseudo-Riemannian geometry encodes information geometry in optimal transport}, Information Geometry. \textbf{5}, 131-159.

\end{thebibliography}
 \end{document}